\documentclass[a4paper,10pt]{article}
\usepackage{amsfonts}
\usepackage{amsmath}
\usepackage{amssymb}
\usepackage{theorem}
\usepackage{amscd}
\usepackage{pstricks}
\usepackage[english]{babel}
\newtheorem{thm}{Theorem}[section]
\newtheorem{lem}[thm]{Lemma}

\newtheorem{prop}[thm]{Proposition}

\newcommand{\res}{\mathrm{Res}}
\newcommand{\ind}{\mathrm{Ind}}
\def\cqfd{\hfill\vrule height 1.6ex width .7ex depth -.1ex }

\def\[{[\![}
\def\]{]\!]}

\title{Modular representations  of cyclotomic Hecke algebras of type $G(r,p,n)$}\author{Gwena\"elle Genet\footnote{Address: {\sc  UFR Math\'ematiques - Case 7012 -      Universit\'e Denis Diderot - Paris 7 -  75251 PARIS Cedex 05, France.}  E-mail address: genet@math.jussieu.fr} and Nicolas Jacon\footnote{
Address:
{\sc Universit\'e de Caen, UFR des Sciences
D\'epartement de Math\'ematiques
Laboratoire LMNO, 
Campus II 14032 Caen cedex, France.}
 E-mail address: jacon@math.unicaen.fr}}

\date{}
\begin{document}
\maketitle
\begin{abstract}
We give a classification of the simple modules for the cyclotomic Hecke algebras over $\mathbb{C}$ in the modular case. We use the unitriangular shape of the
decomposition matrices of Ariki-Koike algebras and Clifford theory.

\end{abstract}
\section{Introduction}

Let $r$, $p$ and $n$ be integers such that $p$ divides $r$. The
complex reflection group of type $G(r,p,n)$ is defined to be the
groups of $n\times n$ permutation matrices such that the entries
are either $0$ or $r^{\textrm{th}}$ roots of unity and the product of the
nonzero entries is a $\frac{r}{p}$-root of unity.

 In  \cite{ariki1} and  \cite{brouemalle}, Ariki and  Brou\'e-Malle have defined
 a Hecke algebra associated  to each complex reflection group
 $G(r,p,n)$. It can be seen as a deformation of the algebra of the group $G(r,p,n)$.
According to a conjecture of Brou\'e and Malle, such algebras,
called cyclotomic Hecke algebras, should occur as endomorphism
algebras of the Lusztig induced character. The representation
theory of cyclotomic Hecke  algebras of type $G(r,1,n)$ (also
known as Ariki-Koike algebras) is beginning to be well understood.
They are  cellular algebras, the simple modules have been
classified in both semi-simple and modular cases and the
decomposition matrices are known in characteristic $0$ (see
\cite{mathas} for a survey of these results).

In \cite{ariki1},  Ariki has shown  that a Hecke algebra of type
$G(r,p,n)$ can be considered as the 0-component of a graded system
for a Hecke algebra of type $G(r,1,n)$ with a special choice of
parameters. As a consequence, in the semi-simple case, he has
given a complete set of non isomorphic simple modules by using
Clifford theory. In the modular case, partial results have been
obtained by Hu in \cite{hu2} (see also  \cite{hu1} and \cite{typedn} where the case
$r=p=2$ which corresponds to Hecke algebras of type $D_n$ is
studied). The main problem is that the restrictions of the simple
modules for non semi-simple Ariki-Koike algebras are  much more
complicated to describe than for semi-simple ones.

The purpose of this paper is precisely to give a parametrization
of the simple modules for the  non semi-simple cyclotomic Hecke algebras
of type $G(r,p,n)$ over $\mathbb{C}$. It comes from the ``canonical basic set''
introduced by Geck and Rouquier in \cite{GR} and extended by the
second author in \cite{papier2}. This set induces a parametrization of
the simple modules of non semi-simple Ariki-Koike algebras by some
FLOTW multipartitions (this kind of multipartitions has been
defined by Foda, Leclerc, Okado, Thibon and Welsh in
\cite{FLOTW}). We proceed in an analogous way than in
\cite[Theorem 2.1]{gwen1}, using the unitriangular shape of the
decomposition matrices of Ariki-Koike algebras and Clifford
theory.

\section{Cyclotomic Hecke algebras, Clifford theory}

\subsection{Cyclotomic Hecke algebras}

Let $r$, $p$ and $n$ be integers such that $p$ divides $r$. We
denote $\displaystyle d:=\frac{r}{p}$. Let $R$ be an integral
ring, $q$ a unit of $R$ and a sequence ${x} = (x_1,\ldots,x_d)$ of elements in $R$.\\
Denote $\mathfrak{H}^{q,x}_{r,p,n} (R)$, the Hecke $R$-algebra
associated to the complex reflection group $G(r,p,n)$, for the
parameters  $q,\ {x}$, defined by the following presentation:\\
$\begin{array}{ll} generators: & a_0,\ldots,a_n,\\
{relations}:&(a_0 - x_1)\ldots(a_0 - x_d) = 0\\
&(a_i -q)(a_i +1) = 0, \ 1 \le i \le n,\\
&a_1 a_3 a_1 = a_3 a_1 a_3, \\
&a_i a_{i+1} a_i = a_{i+1} a_i a_{i+1},\ 2 \le i \le n-1,\\
&(a_1 a_2 a_3)^2 = (a_3 a_1 a_2)^2,\\
&a_1 a_i = a_i a_1,\ 4 \le i \le n,\\
&a_i a_j = a_j a_i,\ 2 \le i < j \le n,\ j \ge i+2,\\
&a_0 a_1 a_2 = (q^{-1} a_1 a_2)^{2-r}a_2 a_0 a_1 + (q-1)
\displaystyle \sum _{i = 1}^{r-2} (q^{-1} a_1 a_2)^{1-k} a_0 a_1,\\
& \hspace{1.3cm} =a_1 a_2 a_0,\\
&a_0 a_i = a_i a_0,\ 3 \le i \le n.\\

\end{array}$\\

We have the following special cases: \begin{itemize} \item if
$r=p=1$, $\mathfrak{H}^{q,x}_{r,p,n} (R)$ is the Hecke algebra of
type $A_{n-1}$, \item if $r=2$ and $p=1$,
$\mathfrak{H}^{q,x}_{r,p,n} (R)$ is the Hecke algebra of type
$B_{n}$, \item if  $r=p=2$,
$\mathfrak{H}^{q,x}_{r,p,n} (R)$ is the Hecke algebra of type
$D_{n}$, \item if $p=1$, $\mathfrak{H}^{q,x}_{r,p,n}(R)$ is  the
Ariki-Koike algebra defined in \cite{arikoi}.\\ \end{itemize}

Suppose the ring $R$ contains a $p^{\textrm{th}}$ root of unity $\eta _p$ and
a $p^{\textrm{th}}$ root $y_i$ of~$x_i$, for each $i \in [1,d]$.  We define a
sequence $Q = (Q_1,\ldots,Q_r)$ of elements in $R$ from the
sequence $x$: for $j \in [1,r]$, $j = sp +t$ with $s \in [0,d-1]$
and $t \in [1,p]$, let  $Q_j := y_{s+1} \eta_p ^{t-1}$. Then,
$\mathfrak{H}^{q,Q}_{r,1,n} (R)$ is the $R$-algebra defined by
$\begin{array}{ll} generators : & T_1,\ldots,T_n ,\\
{relations}: & (T_1 ^p - x_1)\ldots(T_1 ^p - x_d) = 0,\\
& (T_i -q)(T_i +1) = 0, \ 2 \le i \le n,\\
& T_1 T_2 T_1 T_2 = T_2 T_1 T_2 T_1, \\
& T_i T_{i+1} T_i = T_{i+1} T_i T_{i+1},\ 2 \le i \le n-1,\\
& T_i T_j = T_j T_i,\ 1 \le i < j \le n,\ j \ge i+2.\\
\end{array}$\\

We use \cite[Proposition 1.6]{ariki1} and identify $a_0$ with $T_1
^p$, $a_1$ with $T_1^{-1}T_2T_1$ and $a_i$ with $T_i$, for $i \in
[2,n]$, such that $\mathfrak{H}^{q,x}_{r,p,n} (R)$ is considered
as a subalgebra of $\mathfrak{H}^{q,Q}_{r,1,n} (R)$. Then by
\cite[         \S             4.1]{gwenthese}, the Ariki-Koike algebra
$\mathfrak{H}^{q,Q}_{r,1,n} (R)$ is graded over the cyclotomic
Hecke algebra $\mathfrak{H}^{q,x}_{r,p,n} (R)$ by a cyclic group
of order $p$:
$$\mathfrak{H}^{q,Q}_{r,1,n} (R) = \bigoplus _{i = 0} ^{p-1} T_1
^i\, \mathfrak{H}^{q,x}_{r,p,n} (R).$$ \\

Put $\mathfrak{H} := \mathfrak{H}^{q,Q}_{r,1,n} (R)$ and
$\mathfrak{H}' := \mathfrak{H}^{q,x}_{r,p,n} (R)$. We can define
an automorphism $g_{\mathfrak{H}'}^{\mathfrak{H}}$ of
$\mathfrak{H}'$ by
$$g_{\mathfrak{H}'}^{\mathfrak{H}}(h):=T_1^{-1} h T_1,\ \mathrm{for\ all\ }h \in
\mathfrak{H}'.$$  We also have an automorphism
$f_{\mathfrak{H}'}^{\mathfrak{H}}$ of $\mathfrak{H}$ which is
defined on the grading by
$$f_{\mathfrak{H}'}^{\mathfrak{H}}({T_1}^j h)=\eta_p^j\, T_1^j h,\ \mathrm{for\ }j
\in [0,p-1],\ h\in \mathfrak{H}'.$$\\

The next theorem of Dipper and Mathas will be very useful in the following of this
paper. \\

\begin{thm}[{\cite{dipmat}}] \label{theodipmat}
With the previous notations, assume that we have a partition of
$Q$:
$$Q = Q^1\coprod{Q^2}\coprod \ldots \coprod{Q^{s}}$$
such that
$$f_{\Gamma}(q,Q)=\prod_{1\leq{\alpha}<\beta\leq{s}}
{\prod_{{Q_i\in{Q^{\alpha}}\atop{Q_j\in{Q^{\beta}}}}}}{\prod_{-n<a<n}{(q^aQ_i-Q_j)}}$$
is a unit of $R$. Then, ${\mathfrak{H}^{q,Q}_{r,1,n} (R)}$ is
Morita-equivalent to the following algebra:
$$\bigoplus_{{n_1,...,n_s>0}\atop{n_1+...+n_s=n}}{\mathfrak{H}^{q,Q^1}_{r_1,1,n_1} (R)}
\otimes \ldots \otimes{\mathfrak{H}^{q,Q^s}_{r_s,1,n_s} (R)}$$
where for $i \in [1,s]$, $r_i = |Q^i|$. \\
\end{thm}

\subsection{Clifford theory} \label{nota}

\indent Let ${\mathbf y}_1$, ${\mathbf y}_2$,...,${\mathbf y}_{d}$
and ${\mathbf v}$ be indeterminates over $\mathbb{C}$. For $i \in
[1,d]$, put
$${\mathbf x}_i:={\mathbf y}_i^p$$ and for $i \in [1,r],$ $i = sp+t$ with
$s \in [0,d-1]$, $t \in [1,p]$, $\eta _p := \exp
(\frac{2i\pi}{p})$,
$${\mathbf Q}_i := \eta _p ^{t-1} {\mathbf y}_{s+1}.$$ Form the sequences
$\mathbf{x} := (\mathbf{x}_1,\ldots,\mathbf{x}_d)$ and $\mathbf{Q}
:=
(\mathbf{Q}_1,\ldots,\mathbf{Q}_r)$.\\

Let $A:=\mathbb{C}[{\mathbf x}_1,{\mathbf
x}_1^{-1},\ldots,{\mathbf x}_{d},{\mathbf x}_{d}^{-1},{\mathbf
v},{\mathbf v}^{-1}]$, $K:= \mathbb{C}({\mathbf
x}_1,\ldots,{\mathbf x}_{d},{\mathbf v})$ its field of fractions. \\
\indent Let $\theta : A \to \mathbb{C}$ be a  ring  homomorphism
such that $\mathbb{C}$ is the field of fractions of $\theta (A)$. \\ We put for $i \in [1,d]$,
$$x_i:=\theta ({\mathbf x_i}),$$  for $i\in [1,r]$, $$Q_i:=\theta
({\mathbf Q_i})$$ and
$$v:=\theta ({\mathbf v}).\\ $$

Note that for $i \in  [1,r]$, $i = sp+t$ with $s \in [0,d-1]$,
$t \in [1,p]$, then ${Q}_i := \eta _p ^{t-1} { y}_{s+1}$,  for a
complex number ${ y}_{s+1}$ such that ${ y}_{s+1}^p = { x}_{s+1}$.
\\ We also form the sequences ${x} := ({x_1},\ldots,{x_d})$ and
${Q} := ({Q_1},\ldots,{Q_r})$. \\

We will note for short $\mathbf {H}$, resp. ${H}$, for the
Ariki-Koike algebras $\mathfrak{H}^{\mathbf{v},\mathbf{Q}}_{r,1,n}
(K)$, resp. $\mathfrak{H}^{{v},{Q}}_{r,1,n} (\mathbb{C})$ and
$\mathbf {H '} $, resp. ${H}'$, for the cyclotomic Hecke algebras
$\mathfrak{H}^{\mathbf{v},\mathbf{x}}_{r,p,n} (K)$, resp.
$\mathfrak{H}^{{v},{x}}_{r,p,n} (\mathbb{C})$. \\

Since $\mathbf {H}$, resp. ${H}$, is graded over $\mathbf {H'}$,
resp. ${H}'$, we will give some results
which comes from Clifford theory.\\

Let $\mathfrak{H}$ be the Ariki-Koike algebra $\mathbf{H}$, resp.
${H}$, and $\mathfrak{H}'$ be the cyclotomic Hecke algebra
$\mathbf{H'}$, resp. ${H}'$. In both cases, we will denote by $\res$
the functor of restriction from the category of
the $\mathfrak{H}$-modules to the category of the $\mathfrak{H}'$-modules
and by $\ind$ the functor $\mathfrak{H} \otimes _{\mathfrak{H}'} -
$ from the category of the $\mathfrak{H}'$-modules to the category
of the
$\mathfrak{H}$-modules. \\

Let $V$ be a $\mathfrak{H}$-simple module and $U$ be a simple
submodule of $\res\, V$. Put $$o_{f,V} :=
\min{\{k\in{\mathbb{N}_{>0}}\ |\ {}^{f^k}\! V\simeq{V}\}}$$ where
$f := f^{\mathfrak{H}}_{\mathfrak{H}'}$ and for $i \in
[0,o_{f,V}-1]$, ${}^{f^i}\!V$ is the $\mathfrak{H}$-module with
the same underlying space as~$V$ and the structure of module is
given by composing the original action of $\mathfrak{H}$
with~$f^i$. \\ In the same way, let
$$o_{g,U} := \min{\{k\in{\mathbb{N}_{>0}}\ |\ {}^{g^k}\! U\simeq{U}\}}$$
where $g := g^{\mathfrak{H}}_{\mathfrak{H}'}$ and for $i \in
[0,o_{g,U}-1]$, ${}^{g^i}\!U$ is the $\mathfrak{H}'$-module with
the same underlying space as~$U$ with
the action of $\mathfrak{H}'$ twisted by $g^i$.\\
Denote $[X : Y]$ the multiplicity of a simple module $Y$ in the
semi-simple module~$X$. \\

\begin{lem} \label{lemsimple} Denote for short $g :=
g^{\mathfrak{H}}_{\mathfrak{H}'}$ and $f :=
f^{\mathfrak{H}}_{\mathfrak{H}'}$. Every simple
$\mathfrak{H}'$-module~$U$ appears as a direct summand of the
restriction $\res\, V$ of a simple $\mathfrak{H}$-module~$V$ with
multiplicity-free. With such modules $U$ and $V$, a conjugate
$^{g^i}\!U$ of $U$, $i \in [0,o_{g,U} -1]$, occurs in the
restriction of the conjugate $^{f^j}\!V$ of $V$, for all $j \in
[0,o_{f,V} -1]$ and the conjugates of $V$ are the only simple
$\mathfrak{H}$-modules whose restrictions contain such a
$^{g^i}\!U$. Besides
\begin{eqnarray}
\res\, V & \simeq & \bigoplus _{i = 0} ^{o_{f,V} -1}
{}^{f^i}\!V,\\
\ind\, U & \simeq & \bigoplus _{i = 0} ^{o_{g,U} -1}
{}^{g^i}\!U,\\
p & = & o_{g,U} o_{f,V}.
\end{eqnarray} \end{lem}

\emph{Proof}: It is easy to see that
\begin{equation} \label{eqresV}
\res\, V \simeq [\res\, V : U] \bigoplus _{i = 0} ^{o_{g,U} -1}
{}^{g^i}\!U,\end{equation} (see \cite[Theorem 11.1]{currein}). Using \cite[Proposition 2.2]{gwen1}, we get
$$\ind\, U \simeq [\ind\, U : V] \bigoplus _{i = 0} ^{o_{f,V} -1}
{}^{f^i}V$$ as well as
\begin{equation} \label{eqoVoU} \frac{p}{o_{f,V} o_{g,U}} =
[\res \, V : U] [\ind \, U : V].
\end{equation}
As $\res\, V$ is $g$-stable, \cite[Proposition 11.14]{currein}
implies that $\mathrm{End}_{\mathfrak{H}} (\ind\, \res\, V )$ is
graded over $\mathrm{End} _{\mathfrak{H}'} (\res \, V )$ by a
cyclic group of order $p$. By (\ref{eqresV}) and (\ref{eqoVoU}),
the dimension of $\mathrm{End}_{\mathfrak{H}} (\ind\, \res\, V )$
over $\mathbb{C}$ is $\frac{p^2}{o_{f,V}}$ and by the grading over
$\mathrm{End} _{\mathfrak{H}'} (\res \, V )$, this dimension is
equal to $p [\res\, V : U]^2 o_{g,U}$. By comparison with
(\ref{eqoVoU}), it follows that \begin{equation} \label{eqegmult}
[\res\ V : U] = [\ind\ U : V].\end{equation} Now, by definition of
$o_{f,V}$, for short $\mathfrak{o}$, there exists an
$\mathfrak{H}$-automorphism $$t : V \to
{}^{f^{\mathfrak{o}}}\!V.$$ So, $t^{\frac{p}{\mathfrak{o}}}$ is an
$\mathfrak{H}$-isomorphism. By Schur's lemma,
$t^{\frac{p}{\mathfrak{o}}}$ is a scalar. There exists a complex
number $c$ such that $(\frac{t}{c})^{\frac{p}{\mathfrak{o}}} =
\mathrm{Id}_V$. Then
$$\res\, V \simeq \bigoplus _{i = 0}^{\frac{p}{\mathfrak{o}}-1}
\mathrm{Ker} (\frac{t}{c} - \eta _{{p}} ^{\mathfrak{o} i}
\mathrm{Id}_V).$$ It is easy to see that for each $i \in
[0,\frac{p}{\mathfrak{o}}-1]$, $\mathrm{Ker} (\frac{t}{c} - \eta
_{{p}} ^{\mathfrak{o} i} \mathrm{Id}_V)$ is an
$\mathfrak{H}'$-module isomorphic to ${}^{g^i} \mathrm{Ker}
(\frac{t}{c} - \mathrm{Id}_V)$. So the semi-simple module $\res\,
V$ decomposes into more than ${\frac{p}{\mathfrak{o}}}$ simple
modules and by (\ref{eqoVoU}) by less than
${\frac{p}{\mathfrak{o}}}$ simple modules. Therefore $\res\, V$
decomposes into exactly ${\frac{p}{\mathfrak{o}}}$ simple modules
and $[\ind\, U : V] = 1$ so $[\res\, V : U] = 1$ by (\ref{eqegmult}). \\
As the module induced from a simple $\mathfrak{H}'$-module to
$\mathfrak{H}$ is semi-simple, the assertions of the Lemma are
proved. \cqfd \\

\subsection{The semi-simple case}

We keep the notations of the previous subsection.\\

We study the representation theory of the  algebras $\mathbf{H} :=
\mathfrak{H}^{\mathbf{v},\mathbf{Q}}_{r,1,n} (K)$ and $\mathbf{H}'
:= \mathfrak{H}^{\mathbf{v},\mathbf{x}}_{r,p,n} (K)$. By
\cite[Main Theorem]{ariki2}, $\mathbf {H}$ is a split semi-simple
algebra. This
implies  that $\mathbf {H'}$ is also  split semi-simple.\\

A complete set of non isomorphic modules for $\mathbf {H}$ has
been given in \cite{arikoi} and~\cite{djmathas}. Let $\Pi^r_n$ be
the set of $r$-partitions of $n$. We say that $\lambda =
\big(\lambda {(i)}\big)_{i \in [1,r]}$ is a $r$-partition of
size $n$ if for every $i$, $\lambda {(i)}$ is a partition and
$\displaystyle \Sigma _{i = 1}^r |\lambda{(i)} | = n$. \\
For each $r$-partition ${\lambda}\in{\Pi_n^r}$, we can associate
an $\mathfrak{H}^{\mathbf{v},\mathbf{Q}}_{r,1,n}(A)$-module
$S^{\mathbf{v},\mathbf{Q}}(\lambda)$ called a ``Specht module''
that is free over $A$. Here, we use the definition of
``classical'' Specht module contrary to \cite{djmathas} where the
dual Specht modules are
defined. We have the following theorem.\\

\begin{thm}[{\cite{arikoi},\cite{djmathas}}]
The following set is a complete set of non isomorphic absolutely
irreducible
$\mathfrak{H}^{\mathbf{v},\mathbf{Q}}_{r,1,n}(K)$-modules:
$$\left\{ S^{\mathbf{v},\mathbf{Q}}_K(\lambda):=K\otimes_A
S^{\mathbf{v},\mathbf{Q}}(\lambda) \ |\
{\lambda}\in{\Pi_n^r}\right\}.$$
\\ \end{thm}

For each ${\lambda}\in{\Pi_n^r}$, the set of standard tableaux of
shape $\lambda$ is a basis of the underlying $\mathbb{C}$-vector
space of $S^{\mathbf{v},\mathbf{Q}}_K(\lambda)$ (see
\cite{arikoi}). Using \cite[Propositions 3.16, 3.17]{arikoi} and
some combinatorial
properties, it is easy to get the \\

\begin{prop} \label{prophlambda} Let ${\lambda}\in{\Pi_n^r}$. The map
\begin{equation} \label{eqnvarpi} h^{\mathbf{v},\mathbf{Q}}_\lambda :
{}^{f^{\mathbf{H}}_{\mathbf{H}'}} \!
S^{\mathbf{v},\mathbf{Q}}_K(\lambda) \to
S^{\mathbf{v},\mathbf{Q}}_K \big( \varpi(\lambda) \big)
\end{equation} that sends a standard $\lambda$-tableau $\mathbb{T}
= (\mathbb{T}_i)_{i \in [1,r]}$ to the standard
$\varpi(\lambda)$-tableau $\varpi(\mathbb{T}) :=
(\mathbb{T}_{\varpi ^{-1} (i)}) _{i \in [1,r]}$ is an
$\mathbf{H}$-isomorphism where $\varpi $ is the permutation of
$\mathfrak{S}_r$ defined by $ \mathbf{Q} _{\varpi (i)} : = \eta _p
\mathbf{Q}_i,$ for all $i \in [1,r]$ and $\varpi (\lambda) =
(\lambda
{(\varpi ^{-1} (1) )}, \ldots , \lambda{(\varpi ^{-1} (r))})$.\\
\end{prop}

Define $$o_\lambda := \min{\{k\in{\mathbb{N}_{>0}}\ |\ \varpi^k
(\lambda) =  \lambda\}}.$$ It is clear that $o_\lambda =
o_{f,S^{\mathbf{v},\mathbf{Q}}_K(\lambda)}$ and $$\res\,
S^{\mathbf{v},\mathbf{Q}}_K(\lambda) \simeq \bigoplus _{i=
0}^{\frac{p}{o_\lambda}-1} \mathrm{Ker}(h_\lambda ^{o_\lambda} -
\eta _p ^{o_\lambda i}
\mathrm{Id}_{S^{\mathbf{v},\mathbf{Q}}_K(\lambda)}),$$ 
where $h_\lambda^{o_\lambda}: 
{}^{f^{o_{\lambda}}} \!
S^{\mathbf{v},\mathbf{Q}}_K(\lambda) \to
S^{\mathbf{v},\mathbf{Q}}_K \big(\lambda \big)$ sends a  standard $\lambda$-tableau $\mathbb{T}$ to the standard $\lambda$-tableau $\varpi^{o_{\lambda}}(\mathbb{T})$.  For $ i \in
[0,{\frac{p}{o_\lambda}}-1]$, put
\begin{equation} \label{eqdefS} S^{\mathbf{v},\mathbf{Q}}_K(\lambda ,i) := \mathrm{Ker}(h_\lambda
^{o_\lambda} - \eta _p ^{o_\lambda i}
\mathrm{Id}_{S^{\mathbf{v},\mathbf{Q}}_K(\lambda)}).\end{equation}

Let $\mathcal{L}$ be a set of representatives of
$r$-partitions for the action of the cyclic group generated
by $\varpi$ over
$\Pi_n^r$. It follows from  Lemma \ref{lemsimple}, the \\

\begin{prop} \label{propensH'} The set
$$\{S^{\mathbf{v},\mathbf{Q}}_K(\lambda ,i)\ ;\ i \in [0,{\frac{p}{o_\lambda}}-1],\
\lambda \in \mathcal{L}\}$$ is a complete set of non isomorphic
simple $\mathbf{H}'$-modules.\\ \end{prop}

If $\mathbb{T}$ is a standard $\lambda$-tableau, then for $ j \in
[0, \frac{p}{o_\lambda} -1]$, the element $$(\mathbb{T},j) := \sum
_{i=0}^{\frac{p}{o_\lambda} -1} \eta _p ^{- o_\lambda ij} \varpi
^{o_\lambda i}(\mathbb{T})$$ belongs to the eigenspace
$\mathrm{Ker}(h_\lambda ^{o_\lambda} - \eta _p ^{o_\lambda j}
\mathrm{Id}_{S^{\mathbf{v},\mathbf{Q}}_K(\lambda)})$. It comes
from the fact that $\mathbb{T} = \displaystyle \frac{o_\lambda}{p}
\sum _{j=0}^{\frac{p}{o_\lambda} -1} (\mathbb{T},j)$, that for
every $ j \in [0, \frac{p}{o_\lambda} -1]$, the eigenspace
$\mathrm{Ker}(h_\lambda ^{o_\lambda} - \eta _p ^{o_\lambda j}
\mathrm{Id}_{S^{\mathbf{v},\mathbf{Q}}_K(\lambda)})$ is generated by the
$(\mathbb{T},j)$, $\mathbb{T}$ a standard $\lambda$-tableau. So we
have just described the simple $\mathbf{H} '$-modules, see also
\cite{ariki1}. \\

The following part is now concerned with the non semi-simple case: simple $\mathfrak{H}^{q,Q}_{r,1,n}(\mathbb{C})$-modules and
$\mathfrak{H}^{q,x}_{r,p,n}(\mathbb{C})$-modules.\\

\subsection{The non semi-simple case}

We keep the notations of \S\ref{nota}. \\

\subsubsection{Decomposition maps}

The aim of this part is to define the decomposition maps for $H :=
\mathfrak{H}^{q,Q}_{r,1,n}(\mathbb{C})$ and $H'
:=\mathfrak{H}^{q,x}_{r,p,n}(\mathbb{C})$.
Recall that $A:=\mathbb{C}[{\mathbf x}_1,{\mathbf
x}_1^{-1},\ldots,{\mathbf x}_{d},{\mathbf x}_{d}^{-1},{\mathbf
v},{\mathbf v}^{-1}]$, $K:= \mathbb{C}({\mathbf
x}_1,\ldots,{\mathbf x}_{d},{\mathbf v})$ and that
$\theta : A \to \mathbb{C}$ is a ring homomorphism.

Using \cite{geckdecompo}, there exists a discrete valuation ring
$\mathcal{O}$ with maximal ideal $I(\mathcal{O})$ such that
$A \subset{\mathcal{O}}$ and $I(\mathcal{O})\cap
A=\textrm{ker}(\theta)$. Let $k_0 := \mathcal{O}/I(\mathcal{O})$
be the residue field of $\mathcal{O}$. Denote by${\ }^{\overline{\
}} : \mathcal{O} \to k_0$ the canonical map. We obtain the
algebras
$\mathfrak{H}^{\overline{\mathbf{q}},\overline{\mathbf{Q}}}_{r,1,n}(k_0)$
and
$\mathfrak{H}^{\overline{\mathbf{q}},\overline{\mathbf{x}}}_{r,p,n}(k_0)$.
The field $k_0$ can be considered as an extension of $\mathbb{C}$
that is the quotient field of $\theta (A)$. By \cite[Lemma
7.3.4]{gecklivre}, there is an isomorphism between the Grothendieck
groups of finitely generated $H$-modules and
$\mathfrak{H}^{\overline{\mathbf{q}},\overline{\mathbf{Q}}}_{r,1,n}(k_0)$-modules
(resp. $H'$-modules  and
$\mathfrak{H}^{\overline{\mathbf{q}},\overline{\mathbf{x}}}_{r,p,n}(k_0)$-modules):
$$R_0(H)\simeq
R_0(\mathfrak{H}^{\overline{\mathbf{q}},\overline{\mathbf{Q}}}_{r,1,n}(k_0))
\qquad{\textrm{and}\qquad{R_0(H')\simeq
R_0(\mathfrak{H}^{\overline{\mathbf{q}},\overline{\mathbf{x}}}_{r,p,n}(k_0))}}.$$
As a consequence, we obtain well-defined decomposition maps  by
choosing $\mathcal{O}$-forms for the simple
$\mathfrak{H}^{{\mathbf{q}},{\mathbf{Q}}}_{r,1,n}(K)$-modules
(resp.
$\mathfrak{H}^{{\mathbf{q}},{\mathbf{x}}}_{r,p,n}(K)$-modules) and
reducing them modulo $I(\mathcal{O})$:
$$d^{r,1,n} : R_0(\mathbf{H}) \to
R_0(H),$$
$$d^{r,p,n} : R_0(\mathbf{H}') \to
R_0(H').\\ $$

For a simple $\mathbf{H}$-module $V$ and $W$ a simple
$\mathbf{H'}$-module, there exist non negative integers
$d^{(1)}_{V,M}$ with $M$ a simple $H$-module and $d^{(p)}_{W,N}$
with $N$ a simple $H'$-module such that:
\begin{align*}
 d^{r,1,n}([V])&=\sum_{M\in{\mathrm{Irr}(H)}} d^{(1)}_{V,M}[M] \\
\textrm{and}\qquad d^{r,p,n}([W])&=\sum_{N\in{\mathrm{Irr}(H')}}
d^{(p)}_{W,N}[N].\\
\end{align*}

By using the same argument as in \cite[Lemma 5.2]{geck2} (see also
\cite[Theorem 7.4.3.c]{gecklivre}), we obtain the following result.
\begin{prop}\label{commutativite} The following diagram commutes:
\begin{equation*}
\begin{CD}
   R_0(\mathbf{H})       @>\res >>     R_0(\mathbf{H'})      \\
@V{d^{r,1,n}}VV          @VV{d^{r,p,n}}V \\
  R_0(H)                @> \res  >>  R_0(H')
\end{CD}
\end{equation*}
 Moreover, for any simple $\mathbf{H}$-module $V$ and  any simple $H$-module
 $M$, we have
$$d^{(1)}_{V,M}=d^{(1)}_{{}^{f ^{\mathbf{H}}_{\mathbf{H}'}}\! V,{}^{f^H_{H'}}\!M},$$
and for any simple $\mathbf{H'}$-module $W$ and any simple
${H'}$-module $N$, we have:
$$d^{(p)}_{W,N}=d^{(p)}_{{}^{g ^{\mathbf{H}}_{\mathbf{H}'}}\!  W,
{}^{g^H_{H'}}\! N}.\\ $$ \end{prop}

In the following, we will see that under some additionnal
hypotheses, the above property leads to some interesting results
about the decomposition map of cyclotomic Hecke algebras of type
$G(r,p,n)$ and about the simple $H'$-modules.\\

\subsubsection{Simple modules of non semi-simple Ariki-Koike algebras}\label{order}

We will work under the following hypothesis.
We assume that $v$ is a primitive root of unity, $$v := \eta_e =
\exp (\frac{2i\pi}{e}),$$ for a positive integer $e
>1$. Let $$f := \gcd (e,p),\ e = f e',\ p = f p',$$ then $\eta _f
= \eta _p ^{p'} = \eta _e ^{e'}$.\\
By the definition of $Q$ and Theorem \ref{theodipmat}, without
lost of generality, we can suppose that there exist integers
$v_1\leq \ldots\leq v_{\delta}$ that belong to $[0,e'-1]$ such that $Q$
consists of the complex numbers $\eta _e ^{v_i} \eta_p ^j$, where
$i$ ranges over $[1,\delta]$ and $j$ ranges over $[0,p-1]$. Then $r =
\delta fp'$ and we can split $Q$ into $p'$ sets:
$$Q = Q^1 \coprod \ldots \coprod Q^{p'},$$ with $Q^j$ formed by
the $\eta _e ^{v_i} \eta_p ^{lp'+j}$, $i \in [1,\delta]$ and $l \in
[0,f-1]$. We consider $Q^j$ as an ordered sequence
\begin{eqnarray} \lefteqn{Q^j = (\eta_p ^{j-1}\eta _e ^{v_1}, \ldots,
\eta_p ^{j-1}\eta _e ^{v_{\delta}}, \eta_p ^{p'+j-1}\eta _e ^{v_1},
\ldots, \eta_p ^{p'+j-1}\eta _e ^{v_{\delta}}, \ldots  }\nonumber\\
& & \null \qquad {}\eta_p ^{(f-1)p'+j-1}\eta _e ^{v_1}, \ldots,
\eta_p^{(f-1)p'+j-1}\eta _e ^{v_{\delta}}). \end{eqnarray}\

\noindent \emph{Remark}: First note that, for $j \in [1,p']$,
$Q^j = \eta _p ^{j-1} Q^1$ and $$Q^1 = (\eta _e ^{v_1},\ldots,\eta
_e ^{v_{\delta}}, \eta _e ^{v_1+e'},\ldots,\eta _e ^{v_{\delta} +
e'},\ldots,\eta _e ^{v_1+(f-1)e'},\ldots,\eta _e ^{v_{\delta} +
(f-1)e'}),
$$ with $v_1\leq
\ldots\leq v_{\delta} \leq v_1 + e' \leq \ldots \leq v_{\delta} + (f-1) e'$. \\
\indent Note also that the quotient of two elements that belong to
$Q^j$, for  $j \in [1,p']$, is a power of $v = \eta_e$ while the
quotient of an element of~$Q^{j}$ by an element of~$Q^l$, $l \in
[1,p']$,
$l \ne j$ is not such a power, by definition of $f$.\\

If we order the elements of $Q$ by the ordered elements of $Q^1$,
then the ordered elements of $Q^2$, etc, Theorem \ref{theodipmat}
implies that the Ariki-Koike algebra $H$ is Morita equivalent to
the algebra
$$\bigoplus _{{n_1,\ldots,n_{p'}\geq 0}\atop{n_1+...+n_{p'}=n}}
\mathfrak{H}^{q,Q^1}_{f\delta,1,n_1} (\mathbb{C}) \otimes \ldots
\otimes \mathfrak{H}^{q,Q^{p'}}_{f\delta,1,n_{p'}}(\mathbb{C}).$$  So
it is clear that $H$ is also Morita equivalent to the algebra
$$\bigoplus _{{n_1,\ldots,n_{p'}\geq 0}\atop{n_1+...+n_{p'}=n}}\
\mathfrak{H}^{q,Q^1}_{f\delta,1,n_1} (\mathbb{C}) \otimes \ldots
\otimes \mathfrak{H}^{q,Q^{1}}_{f\delta,1,n_{p'}}(\mathbb{C}).$$\\
We order $\mathbf{Q}$  and we split it into $p'$ sets following
what we have done for $Q$: $$ \mathbf{Q} = \mathbf{Q}^1\coprod
\ldots \coprod \mathbf{Q}^{p'}.\\ $$ \\

\subsubsection{Parametrization of simple $H$-modules}\label{avalue}

We can now give a parametrization for the simple $H$-modules. By
\cite[Theorem~2.5]{arikicyc}, the simple modules of the
Ariki-Koike algebra $\mathfrak{H}^{q,Q_1}_{f\delta,1,n_i}
(\mathbb{C})$, $i \in [1,p']$, are the quotient modules $D^{v,Q^1}
_\mathbb{C}(\lambda) := S^{\mathbf{v},\mathbf{Q^1}} _\mathbb{C}
(\lambda) / J(S^{\mathbf{v},\mathbf{Q^1}} _\mathbb{C} (\lambda) )
$ of the specialisations of the Specht modules
$S^{\mathbf{v},\mathbf{Q^1}} _\mathbb{C} (\lambda) := \mathbb{C}
\otimes S^{\mathbf{v},\mathbf{Q^1}} (\lambda)$ by their radical of
Jacobson, with $\lambda$ a Kleshchev $f\delta$-partition  (for
the parameters $Q^1$) of size $n_i$. So the simple modules of $H$
are labelled by the $p'$-tuples of Kleshchev $f\delta$-partitions
(associated to $Q^1$) $(\lambda ^1,\ldots,\lambda
^{p'})$ with $\displaystyle \sum _{i = 1}^{p'} |\lambda ^i| = n$.
We will denote this set by $\Lambda ^0$ and for ${\lambda} \in
\Lambda ^0$, $ D^{{v},{Q}}
_\mathbb{C} (\lambda)$ the corresponding simple $H$-module. \\

In \cite{papier2}, another parametrization for the simple
$\mathfrak{H}^{q,Q_1}_{f\delta,1,n_i}(\mathbb{C})$-modules
has been found by using Lusztig's $a$-function (see Theorem \ref{ar} below).
Following this paper, we will  associate to each  $r$-partition an $a$-value.
To do this, we first describe the $a$-value of a
$f\delta$-partition $\mu=(\mu(1),\mu(2),...,\mu(f\delta))$.

For $k=1,...,\delta$, $s=1,...,f$ and  $j=(s-1)\delta+k$  we put $w_j:=v_k+(s-1)e'$
so that  $w_j$ is the power of the $j^{\textrm{th}}$-component of ${Q^1}$.

Then, for  $j:=1,...,f\delta$, we define:
$$m^{(j)}:=w_j-\frac{je}{f\delta}+e$$
and for $s=1,...,n$, we put:
$$B'^{(j)}_s:=\mu(j)_s-s+n+m^{(j)}$$
where we use the convention that $\mu(j)_p:=0$ if $p$ is greater than the height
of  $\mu (j)$. For $j=1,...,f\delta$, let $B'^{(j)}=(B'^{(j)}_1,...,B'^{(j)}_n)$.
Then, we define:
$$a^{(1)}_{f\delta}({\mu}):=\sum_{{0\leq{i}\leq{j}<f\delta}\atop{{(a,b)\in{B'^{(i)}
\times{B'^{(j)}}}}\atop{a>b\ \textrm{if}\ i=j}}}{\min{\{a,b\}}}   -\sum_{{1\leq{i,j}
\leq f\delta}\atop{{a\in{B'^{(i)}}}\atop{1\leq{k}\leq{a}}}}{\min{\{k,m^{(j)}\}}}+g(n).$$
where $g(n)$ is a rational number which only depends on the $m^{(j)}$ and on $n$
(the expression of $g$ is given in \cite{papier2}).

 Let  ${\lambda} \in \Pi _n ^r$,
write $\lambda = (\lambda(1),\ldots,\lambda(r))$ as $\lambda =
(\lambda ^1,\ldots,\lambda ^{p'})$ where, for $i \in [1,p']$,
$\lambda ^i = \big(
\lambda(f\delta(i-1)+1),\ldots,\lambda(f\delta(i-1)+f\delta)\big) $.

Finally, we define
$$a^{(1)}_r(\lambda ) := \displaystyle \sum _{i=1}^{p'}
a^{(1)}_{f\delta}(\lambda ^i).$$

To each Specht module $S^{\mathbf{v},\mathbf{Q^1}}_K (\lambda)$ of
$\mathfrak{H}^{\mathbf{q},\mathbf{Q_1}}_{f\delta,1,n'} (K)$, $n' \in
[0,n]$ and $\lambda$ an $f\delta$-partition, we also associate an
$a$-value
$$a^{(1)}_{f\delta}\big(S_K ^{\mathbf{v},\mathbf{Q^1}} (\lambda)\big)
:= a^{(1)}_{f\delta}(\lambda).$$ \\

To each simple $\mathfrak{H}^{{q},{Q_1}}_{f\delta,1,n'}
(\mathbb{C})$-module $D^{{v},{Q^1}} _\mathbb{C} (\lambda)$, $n'
\in [0,n]$ and $\lambda$ a Kleshchev $f\delta$-partition of size
$n'$, we again associate  an $a$-value
$$a_{f\delta}^{(1)}({D^{{v},{Q^1}} _\mathbb{C} (\lambda)}) := \min \{a^{(1)}_{f\delta}(\mu)\ | \ \mu
\in \Pi _{n'}^{f\delta},\ d^{(1)} _{S^{\mathbf{v},\mathbf{Q^1} }_K
(\mu),{D^{{v},{Q^1}} _\mathbb{C} (\lambda)}} \ne 0\}.$$ \\

\begin{thm}[{{\cite[Theorem 2.3.8]{these}}}] \label{ar}Let $\lambda$ be a Kleshchev
$f\delta$-partition of $n'$ for a non-negative integer $n'$. There
exists a unique $\mu := \kappa (\lambda)$ such that
$a_{f\delta}^{(1)}({D^{{v},{Q^1}} _\mathbb{C} (\lambda)}) =
a_{f\delta}^{(1)}(\mu)$ and $d^{(1)}
_{S^{\mathbf{v},\mathbf{Q^1}}_K (\mu),{D^{{v},{Q^1}} _\mathbb{C}
(\lambda)}} \ne 0$. The function $\kappa$ is a bijection between
the set of the Kleshchev $f\delta$-partitions of size $n'$ and the
set of the FLOTW $f\delta$-partitions of size $n'$.\\ \end{thm}

Recall what is a FLOTW $f\delta$-partition $\lambda = \big(\lambda
(1),\ldots,\lambda (f\delta)\big)$ of size~$n'$ associated to the
parameters $v$ and $Q^1$, see \cite{FLOTW}. Denote each partition
$\lambda (i),$ $i \in [1,f\delta]$, as $(\lambda(i)_1,\lambda(i)_2,\ldots)$. The
multipartition $\lambda$ satisfies
\begin{enumerate} \item for $i \in [1,f]$, $j \in [1,\delta-1]$, $k$
positive integer, $$\lambda\big((i-1)\delta+j\big)_k \geq \lambda
\big((i-1)m+j+1\big) _{k + v_{j+1} - v_j},$$ \item for $i \in
[1,f-1]$, $k$ positive integer, $$\lambda (i\delta)_k \geq \lambda
(i\delta+1)_{k + v_{1} + e' - v_{\delta}},$$ \item   for $k$ positive integer,
$$\lambda (f\delta)_k \geq \lambda (1) _{k + v_{1} + e' - v_{\delta}},$$ as
well as \item to each node of the diagram of $\lambda$ which is located in
the $a^{\textrm{th}}$ row and  the $b^{\textrm{th}}$ column of the $c^{\textrm{th}}$
partition of $\lambda$ (for two poditive integers $a$ and $b$ and $1\leq c\leq f\delta$),
we associate its residue $\eta_e^{b-a+v_c}$. Then, for any positive integer $k$,
the cardinality of the
set of the residues associated to the nodes in both the $k^{\textrm{th}}$
columns and the right rims of the Young diagrams of $\lambda$ is
not $e$.
\\\end{enumerate}

In the same way, to each Specht module $S^{\mathbf{v},\mathbf{Q}}
_K (\lambda)$ of $\mathbf{H}$, $\lambda \in \Pi _n ^r$, $\lambda =
(\lambda ^1,\ldots,\lambda ^{p'})$, we can define
$a_r^{(1)}(S^{\mathbf{v},\mathbf{Q}} _K (\lambda)) :=
a_r^{(1)}(\lambda)$ that
is the sum $\displaystyle \sum _{i =1}^{p'} a^{(1)}_{f\delta}(\lambda^i)$.\\
To each simple $H$-module $D^{v,Q} _{\mathbb{C}}(\lambda )$,
$\lambda \in \Lambda^0$, we associate the $a$-value
\begin{equation} \label{eqndefamu}
a_r^{(1)}({D^{v,Q} _{\mathbb{C}}(\lambda )}) := \min \{a_r^{(1)}(\mu)\
| \ \mu \in \Pi _{n}^{r},\ d^{(1)} _{S^{\mathbf{v},\mathbf{Q}}_K
(\mu),{D^{v,Q} _{\mathbb{C}}(\lambda )}} \ne 0\}.\\
\end{equation}
With Theorem \ref{theodipmat}, it is clear that there exists a
unique $\mu := \kappa ' (\lambda)$ such that $a_r^{(1)}({D^{v,Q}
_{\mathbb{C}}(\lambda )})= a_r^{(1)}(\mu)$ and $d^{(1)}
_{S^{\mathbf{v},\mathbf{Q}}_K (\mu),{D^{v,Q} _{\mathbb{C}}(\lambda
)}} \ne 0$. In fact, $$\kappa ' (\lambda) = (\kappa
(\lambda^1),\ldots,\kappa (\lambda^{p'}))$$ if $\lambda = (\lambda
^1,\ldots,\lambda ^{p'})$.
\\ The function $\kappa '$ is a bijection between $\Lambda ^0$ and
$\Lambda ^1$ the set of the $p'$-tuples of FLOTW
$f\delta$-partitions $(\lambda ^1,\ldots,\lambda ^{p'})$ of size
$\displaystyle \sum _{i=1}^{p'} |\lambda ^i| = n$.\\

\subsubsection{Preliminary results}

In order to prove the main Theorem, we will need some
preliminary results. First, recall that we have ordered ${\bf Q}$ by the ordered
elements of ${\bf Q}^1$, then by the ordered elements of  ${\bf Q}^2$ etc as in
\S \ref{order}. Let $\varpi$ be the permutation of $\mathfrak{S}_r$ defined in
Proposition   \ref{prophlambda} and  let
$$\lambda:=(\lambda^1[1],...,\lambda^1 [f\delta],\lambda^2 [1],...,\lambda^2
[f\delta],...,\lambda^{p'}[1],...,\lambda^{p'} [f\delta])$$
be an $r$-partition. Hence, following the notations of the previous paragraph,
we have $\lambda^j[i]=\lambda(f\delta (j-1)+i)$. Then,  it is easy to verify
that \begin{itemize} \item
for $i \in [1,(p'-1)f\delta]$, then $\varpi(i) = i + f\delta$, \item for $i
\in [r-f\delta +1, r-\delta]$, then $\varpi (i) = i-r+(f+1)\delta$, \item for $i
\in [r-\delta+1, r]$, then $\varpi (i) = i - r + \delta$.
\end{itemize}
Then, by Proposition  \ref{prophlambda}, we have:
$$S^{\mathbf{v},\mathbf{Q}}_K(\lambda) \simeq
S^{\mathbf{v},\mathbf{Q}}_K \big( \varpi(\lambda) \big)$$
where $\varpi(\lambda)$ $=(\lambda^{p'}[f\delta-\delta+1]$ $,\ldots,$ $\lambda^{p'}
[f\delta],$$ \lambda^{p'}[1]$ $,\ldots,$ $\lambda^{p'}[f\delta-\delta],$ $\lambda^1[1]$
$,\ldots,$ $\lambda^1 [f\delta],$ $\ldots$ $,\lambda^{p'-1}[1]$ $,\ldots,$
$\lambda^{p'-1} [f\delta])$.

\begin{prop} \label{propegaSpechtf} Let $\lambda \in \Pi _n ^r$,
then $a^{(1)}_r(S^{\mathbf{v},\mathbf{Q}}_K (\lambda)) =
a^{(1)}_r({}^{f^{\mathbf{H}}_{\mathbf{H}'}}\! S^{\mathbf{v},\mathbf{Q}}_K
(\lambda))$. Besides, if $\lambda \in \Lambda ^0$, then
$a_r^{(1)}({D^{v,Q} _{\mathbb{C}} (\lambda)}) =
a_r^{(1)}({{}^{f^H_{H'}}\! D^{v,Q} _{\mathbb{C}} (\lambda)})$.\\
\end{prop}

\emph{Proof}:
By the definition of the $a$-value (see \S \ref{avalue}), we have:
$$a^{(1)}_r(S^{\mathbf{v},\mathbf{Q}}_K (\lambda)):=\sum_{i=1}^{p'} a^{(1)}_{f\delta}
(\lambda^i[1],...,\lambda^i [f\delta])$$
Thus, by using the above properties, it is sufficient to show that $a^{(1)}_{f\delta}$
$(\lambda^{p'}[1]$ $,\ldots$ $,\lambda^{p'} [f\delta])$$=a^{(1)}_{f\delta}( \lambda^{p'}
[f\delta-\delta+1],\ldots,\lambda^{p'} [f\delta],\lambda^{p'}[1],\ldots,\lambda^{p'}
[f\delta-\delta])$.

Using the notations of \S \ref{avalue}, for $i=1,...,\delta$ and $j=1,...,f$
we write $m^{(i)}[j]:=m^{((j-1)\delta+i)}$. Then, all we have to do is to prove
that $m^{(i)}[j]$ doesn't depend on $j$. We have  $\displaystyle m^{(i)}[j]=v_i-\frac{((j-1)
\delta+i)e}{f\delta}+(j-1)e'=v_i-\frac{ie}{f}$ because $e'f=e$.
Hence, we obtain the desired result.
\\

With the definition (\ref{eqndefamu}), it is trivial that the
previous result about Specht modules implies that, for any
$\lambda \in \Lambda ^0$, $a_{r}^{(1)}({D^{v,Q} _{\mathbb{C}}
(\lambda)}) = a_r^{(1)}({{}^{f^H_{H'}}\! D^{v,Q}
_{\mathbb{C}} (\lambda)})$. \cqfd \\

With Proposition \ref{propegaSpechtf} and Lemma \ref{lemsimple},
we can associate an $a$-value to each simple $\mathbf{H'}$-module
$U$. If $U$ appears in the restriction $\res\,
S^{\mathbf{v},\mathbf{Q}}_K (\lambda)$, for $\lambda \in \Pi _n
^r$, put $$a_r^{(p)} (U) := a_r^{(1)}\big(S^{\mathbf{v},\mathbf{Q}}_K
(\lambda)\big). $$\\

 In the same way, we can associate an $a$-value to
each simple ${H'}$-module $W$. If $W$ appears in the restriction
$\res\, D^{v,Q} _{\mathbb{C}}(\lambda)$, for $\lambda \in \Lambda
^0$, then
$$a_r^{(p)} (W) := a_r^{(1)}({D^{v,Q}
_{\mathbb{C}}(\lambda)}).$$ \\

Another useful proposition is the following one.\\

\begin{prop} \label{propstab}
Let $\lambda \in \Lambda ^1$, then $\varpi(\lambda) \in \Lambda
^1$, where $\varpi$ is defined in the \emph{Proposition
\ref{prophlambda}}. \\ \end{prop}

\emph{Proof}:
Write $\lambda \in \Lambda ^1$ as $\lambda = (\lambda
^1,\ldots,\lambda ^{p'})$ where $\lambda ^i = (\lambda ^i
[1],\ldots,\lambda ^i [f\delta])$ is a FLOTW $f\delta$-partition for
$Q^1$. We will see that $\varpi (\lambda) \in \Lambda ^1$ so
the result will be proved. The $r$-partition $\varpi
(\lambda)$ is equals to  $(\lambda^{p'}[f\delta-\delta+1]$ $,\ldots,$
$\lambda^{p'}[f\delta],$$ \lambda^{p'}[1]$ $,\ldots,$ $\lambda^{p'}[f\delta-\delta],$
$\lambda^1[1]$ $,\ldots,$ $\lambda^1 [f\delta],$ $\ldots$ $\lambda^{p'-1}[1]$
$,\ldots,$ $\lambda^{p'-1} [f\delta])$. It belongs to $\Lambda ^1$ if and only if
$(\lambda^{p'}[f\delta-\delta+1]$ $,\ldots,$ $\lambda^{p'}[f\delta],$
$ \lambda^{p'}[1]$ $,\ldots,$ $\lambda^{p'}[f\delta-\delta]$)
is a FLOTW $f\delta$-partition. It
is easy to verify that the three first conditions hold. Now, the
set of the residues associated to the nodes in both the $k^{\textrm{th}}$
columns and the right rims of the Young diagrams of this last
multipartition consists of the residues of $\lambda^{p'}$ but
multiplied by $\eta ^{e'} _e$. So the fourth condition about
FLOTW multipartitions hold as well as the Proposition. \cqfd \\

\section{A parametrization of the simple  $H'$-modules  }

We keep the above notations. The aim of this section is to
describe a parametrization of the simple $H'$-modules, with $H' =
\mathfrak{H}^{v,Q}_{r,p,n} (\mathbb{C})$. The proof is highly
inspired by \cite[Theorem 2.1]{gwen1}.\\

\begin{thm}  \label{theofinal} For all simple $H'$-module $N$,
there exists a simple $\mathbf{H'}$-module $W_N$ such that
$$d^{r,p,n}([W_N])=[N]+\sum_{{a_r^{(p)}(L)<a_r^{(p)}(N)}}d_{{W_N},L}[L],$$
and $a_r^{(p)}(N) = a_r^{(p)}(W_N)$.\\
Let ${\Lambda^1}' := \{W_N\ |\ N\in \mathrm{Irr}(H ')\}$. Then
${\Lambda^1}' $ consists in the
$S^{\mathbf{v},\mathbf{Q}}_K(\lambda , i)$, with $\lambda \in
\Lambda ^1 \cap \mathcal{L}$, $i \in [0,\frac{p}{o_\lambda} -1]$
with the notations of the Proposition \ref{propensH'} and the map
$ N\in \mathrm{Irr}(H') \mapsto{W_N \in {\Lambda^1}'}$ is
bijective. So ${\Lambda^1}'$ is a
parametrization of the simple $H'$-modules. \\
\end{thm}

\emph{Proof }: Fix a simple $H'$-module $N$. By Lemma
\ref{lemsimple}, there exists $\lambda \in \Lambda ^0$ such that
$N$ appears in the restriction $\res D^{v,Q} _\mathbb{C} (\lambda
)$. We have \begin{equation} \label{eqd1}
d^{r,1,n}([S^{\mathbf{v},\mathbf{Q}}_{\mathbb{C}}\big(\kappa
(\lambda)\big)]) = [ D^{v,Q} _{\mathbb{C}}(\lambda)] + \sum
_{{a_r^{(p)}({D^{v,Q} _{\mathbb{C}}(\mu)}) <}\atop{ a_r^{(p)}({D^{v,Q}
_{\mathbb{C}}(\lambda)}})}
d^{(1)}_{S^{\mathbf{v},\mathbf{Q}}_{K}(\kappa (\lambda)), D^{v,Q}
_{\mathbb{C}}(\mu)} [D^{v,Q} _{\mathbb{C}}(\mu)],\end{equation}
and $a^{(1)}_r\big(S^{\mathbf{v},\mathbf{Q}}_{\mathbb{C}}(\kappa
(\lambda))\big) = a_r^{(1)}({D^{v,Q} _{\mathbb{C}}(\lambda)})$.  With
the notations of (\ref{eqdefS}) and
Proposition~\ref{commutativite}, this equality implies that
\begin{equation} \label{eqd2}
\sum _{i = 0} ^{\frac{p}{o_{\kappa(\lambda)}}-1}
d^{r,p,n}([S^{\mathbf{v},\mathbf{Q}} _\mathbb{C} (\kappa
(\lambda),i)]) = \sum _{i = 0} ^{\frac{p}{o_{f^H_{H'},D^{v,Q}
_{\mathbb{C}}(\lambda)}}-1}[ {}^{(g^H_{H'}) ^i}\! N] +
[N'],\end{equation} where $N'$ is a sum of $H'$-modules, $D$, with
$a$-value $a_r^{(p)}(D) < a_r^{(p)}({D^{v,Q} _{\mathbb{C}}(\lambda)})$
and, for all~$i$,~$j$, $a_r^{(p)}(S^{\mathbf{v},\mathbf{Q}}
_\mathbb{C} (\kappa (\lambda),i)) = a_r^{(p)}({{}^{(g^H_{H'})
^i}\! N})$.\\ Denote by $\tau(\lambda) \in \Lambda ^0$, the
multipartition defined by $D^{v,Q}_{\mathbb{C}}(\tau(\lambda))
 := {}^{f^H_{H'}}\! D^{v,Q}_{\mathbb{C}}(\lambda) $. If we twist the
action of $H$ by $f^H_{H'}$, the equality (\ref{eqd1}) gives that
$\kappa(\tau(\lambda)) = \varpi (\kappa (\lambda))$ by Proposition
\ref{propstab} and by definition of $\kappa$. As this map is
bijective, it is clear that
$$o_{f^{\mathbf{H}}_{\mathbf{H}'},S^{\mathbf{v},\mathbf{Q}} _K
(\kappa (\lambda))} = o_{f^H_{H'},D^{v,Q}_{\mathbb{C}}(\lambda)
}$$ and therefore for all $i \in
[0,\frac{p}{o_{\kappa(\lambda)}}-1]$,
$$o_{g^{\mathbf{H}}_{\mathbf{H}'},S^{\mathbf{v},\mathbf{Q}} _K
(\kappa (\lambda),i)} = o_{g^H_{H'},N}$$ with the last equality of
Lemma \ref{lemsimple}. \\ Recall that for $i \in
[0,\frac{p}{o_{\kappa(\lambda)}}-1]$, $S^{\mathbf{v},\mathbf{Q}}
_K (\kappa (\lambda),i)$ is conjugate to
$S^{\mathbf{v},\mathbf{Q}} _K (\kappa (\lambda),0)$ by
$(g^{\mathbf{H}}_{\mathbf{H}'})^i$. The equality (\ref{eqd2})
implies that there exists an unique $ i \in
[0,\frac{p}{o_{\kappa(\lambda)}}-1]$ such that
\begin{equation} \label{eqd3} d^{r,p,n}([S^{\mathbf{v},\mathbf{Q}}
_\mathbb{C} (\kappa (\lambda),i)]) = [N] + [N''],\end{equation}
where $N''$ is a sum of $H'$-modules, $D$, with $a$-value
$a_r^{(p)}(D) < a_r^{(p)}({N})$. This construction makes clear all the
statements
of the Theorem. \cqfd \\

\noindent \emph{Remark}: The above theorem is only concerns with
the case where $v$ is a primitive $e^{\textrm{th}}$-root of unity.
If $v$ isn't a root of unity, it is readily checked that an
analogue of this theorem holds by replacing $\Lambda^1$ by the set
of Kleshchev multipartitions $\Lambda^0$ at $e=\infty$. The proof
of this result can be easily obtained following the outline of
\cite{gwen1}.

\end{document}